\documentclass[12pt]{article}  
\usepackage{e-jc}

\usepackage{amsmath} 
\usepackage{latexsym}

\usepackage{url}
\usepackage{amssymb} 
\usepackage{exscale} 
\usepackage{graphicx} 
\usepackage{makeidx}

\usepackage[small,bf,hang]{caption} 

\usepackage{algorithm}
\usepackage{algorithmic}

\newtheorem{thm}{Theorem}

\newtheorem{lem}[thm]{Lemma}

\newenvironment{proof}[1][Proof]
{\par\noindent{\bf #1.} }{\hspace*{\fill}\nolinebreak{$\Box$}\bigskip\par}

\usepackage{color}

\title{\bf New Computational Upper Bounds\\
for Ramsey Numbers $R(3,k)$}

\author{
Jan Goedgebeur\\
\small Department of Applied Mathematics and Computer Science\small \\[-0.8ex]
\small Ghent University, B-9000 Ghent, Belgium\\[-0.8ex]
\small \texttt{jan.goedgebeur@ugent.be}\\
\\
Stanis\l aw P. Radziszowski \\
\small Department of Computer Science \\[-0.8ex]
\small Rochester Institute of Technology, Rochester, NY  14623, USA \\[-0.8ex]
\small \texttt{spr@cs.rit.edu}\\\
}

\date{\dateline{XX}{XX}}

\begin{document}
\maketitle
\begin{abstract}
\noindent
Using computational techniques we derive six new upper
bounds on the classical two-color Ramsey numbers:
$R(3,10) \le 42$,
$R(3,11) \le 50$, $R(3,13) \le 68$, $R(3,14) \le 77$,
$R(3,15) \le 87$, and $R(3,16) \le 98$. All of them are
improvements by one over the previously best published bounds.

Let $e(3,k,n)$ denote the minimum number of edges in any
triangle-free graph on $n$ vertices without independent
sets of order $k$. The new upper bounds on $R(3,k)$ are
obtained by completing the computation of the exact values
of $e(3,k,n)$ for all $n$ with $k \leq 9$ and for
all $n \leq 33$ for $k = 10$,
and by establishing new lower bounds on $e(3,k,n)$
for most of the open cases for $10 \le k \le 15$.
The enumeration of all graphs witnessing the values of $e(3,k,n)$
is completed for all cases with $k \le 9$.
We prove that the known critical graph for $R(3,9)$ on
35 vertices is unique up to isomorphism. For the case of
$R(3,10)$, first we establish that $R(3,10)=43$
if and only if $e(3,10,42)=189$, or equivalently, that if
$R(3,10)=43$ then every critical graph is regular of degree 9.
Then, using computations, we disprove the existence of
the latter, and thus show that $R(3,10) \le 42$.

\end{abstract}

\eject
\section{Definitions and Preliminaries}

\bigskip
In this paper all graphs are simple and undirected. Let $G$ be such a graph.
The vertex set of $G$ is denoted by $V(G)$, the edge set of $G$
by $E(G)$, and the number of edges in $G$ by $e(G)$. The set of
neighbors of $v$ in $G$ will be written as $N_v(G)$ (or just $N(v)$
if $G$ is fixed). The independence
number of $G$, denoted $\alpha(G)$, is the order of the largest
independent set in $G$, $\deg_G(v)$ is the degree of vertex
$v \in V(G)$, and $\delta(G)$ and $\Delta(G)$ are the minimum
and maximum degree of vertices in $G$, respectively.
For graphs $G$ and $H$, $G \cong H$ means that they are isomorphic.

\medskip
For positive integers $k$ and $l$, the \emph{Ramsey number} $R(k,l)$
is the smallest integer $n$ such that if we arbitrarily color
the edges of the complete graph $K_n$ with 2 colors,
then it contains a monochromatic $K_k$ in the first color or
a monochromatic $K_l$ in the second color. If the edges in the
first color are interpreted as a graph $G$ and those in the second color
as its complement $\overline{G}$, then $R(k,l)$ can be defined
equivalently as the smallest $n$ such that every graph on $n$
vertices contains $K_k$ or has independence $\alpha(G)\ge l$.
A regularly updated dynamic survey by the second author \cite{SRN}
lists the values and the best known bounds
on various types of Ramsey numbers.

Any $K_k$-free graph $G$ on $n$ vertices with $\alpha(G)<l$
and $e(G)=e$ will be called a $(k,l;n,e)$-{\em{graph}}, and
by $\mathcal{R}(k,l;n,e)$ we will denote
the set of all $(k,l;n,e)$-graphs.
We will often omit the parameter $e$, or both $e$ and $n$,
or give some range to either of these parameters, when referring
to special $(k,l;n,e)$-graphs or sets $\mathcal{R}(k,l;n,e)$.
For example, a $(k,l)$-graph is a $(k,l;n,e)$-graph for some
$n$ and $e$, and the set $\mathcal{R}(3,9;35,\le 139)$ consists
of all 35-vertex triangle-free graphs with $\alpha(G) \le 8$ and
at most 139 edges (later we will prove that this set is empty).
Any $(k,l;R(k,l)-1)$-graph will be called {\em critical} for $(k,l)$.

\medskip
Let $e(k,l,n)$ denote the minimum number of edges in any
$(k,l;n)$-graph (or $\infty$ if no such graph exists). The sum
of the degrees of all neighbors of $v$ in $G$ will be denoted
by $Z_G(v)$ (or $Z(v)$ if $G$ is fixed), i.e.
$$Z(v)=Z_G(v)=\sum_{\{u,v\}\in E(G)}{\deg_G(u)}.\eqno{(1)}$$

\smallskip
In the remainder of this paper we will study only
triangle-free graphs. Note that for any $G \in \mathcal{R}(3,k)$
we have $\Delta(G)<k$, since all neighborhoods of vertices in $G$
are independent sets.

Let $G$ be a $(3,k;n,e)$-graph.
For any vertex $v \in V(G)$, we will denote by $G_v$
the graph induced in $G$ by the set
$V(G) \setminus (N_G(v) \cup \{v\})$.
If $d=\deg_G(v)$, then
clearly  $G_v$ is a $(3,k-1;n-d-1,e(G)-Z_G(v))$-graph.
Note that this implies that
$$\gamma(v)=\gamma(v,k,G)=e-Z_G(v)-e(3,k-1,n-d-1)\ge 0,\eqno{(2)}$$
where $\gamma(v)$ is the so called {\em deficiency}
of vertex $v$ \cite{GY}. Finally, the deficiency of
the graph $G$ is defined as
$$\gamma(G)=\sum_{v \in V(G)}{\gamma(v,k,G)} \ge 0.\eqno{(3)}$$

\noindent
The condition that $\gamma(G)\ge 0$ will be often sufficient
to derive good lower bounds on $e(k,l,n)$, though a stronger
condition that all summands $\gamma(v,k,G)$
of (3) are non-negative sometimes
implies even better bounds.
It is easy to compute $\gamma(G)$ just from the degree
sequence of $G$ \cite{GY,GR}. If a $(3,k;n,e)$-graph $G$
has $n_i$ vertices of degree $i$, then
$$\gamma(G)=ne-
\sum_{i}{n_i\big( i^2+e(3,k-1,n-i-1)\big)} \ge 0,\eqno{(4)}$$
\noindent
where $n=\sum_{i=0}^{k-1}{n_i}$ and $2e=\sum_{i=0}^{k-1}{i n_i}$.

\bigskip
\bigskip
\section{Summary of Prior and New Results}

\bigskip
In 1995, Kim \cite{Kim} obtained a breakthrough result by
establishing the exact asymptotics of $R(3,k)$ using
probabilistic arguments. Recently, the fascinating story of
developments and results related to the infinite aspects
of $R(3,k)$ was written by Spencer \cite{Spe}.

\begin{thm}[\cite{Kim}]
\label{theorem:kim}
$R(3,k) = \Theta(n^2/\log n)$.
\end{thm}

\bigskip
Theorem~\ref{theorem:kim} gives the exact asymptotics
of $R(3,k)$, while computing
the values for concrete cases remains an open problem for all
$k \ge 10$. Still, the progress obtained in the last 50 years
in this area is remarkable. Known exact values of $R(3,k)$
for $k \le 9$,
and the best lower and upper bounds for higher $k$, are listed
in \cite{SRN} together with all the references. We note that
much of this progress was obtained with the use of
knowledge about $e(3,k,n)$. This direction is also the main 
focus of our paper: we compute new exact values of $e(3,k,n)$
in several cases and give improved lower bounds for many other,
which in turn permits us to prove new upper bounds on $R(3,k)$
for $k=10, 11, 13, 14, 15$ and $16$. Likely, more new upper bounds
could be obtained for some $17 \le k \le 21$, but we did
not perform these computations.

\bigskip
General formulas for $e(3,k,n)$ are known for all $n \le 13k/4-1$
and for $n=13k/4$ when $k=0 \bmod 4$.

\medskip
\begin{thm}[\cite{RK1,RK3}]
\label{theorem:comulative_small}
For all $n,k \ge 1$, for which $e(3,k+1,n)$ is finite,
$$
e(3,k+1,n) = 
\left \{
\begin{array}{ll}
0 & \textrm{  if  }\  n \le k,\\
n - k & \textrm{  if  }\  k < n \le 2k,\\
3n - 5k & \textrm{  if  }\  2k < n \le 5k/2,\\
5n - 10k & \textrm{  if  }\  5k/2 < n \le 3k,\\
6n - 13k & \textrm{  if  }\  3k < n \le 13k/4-1.
\end{array}
\right. \eqno{(5)}
$$

\noindent
Furthermore, $e(3,k+1,n) = 6n-13k$ for $k=4t$ and $n=13t$,
and the inequality $e(3,k+1,n) \ge 6n-13k$ holds for all
$n$ and $k$. All the critical graphs have been
characterized whenever the equality in the theorem holds
for $n \le 3k$.
\end{thm}

\medskip
Theorem~\ref{theorem:comulative_small} is a cumulative summary of various contributions
\cite{GY, GR, RK1, RK2, RK3}. It captures many of the small cases, as
presented in Table 3 in Section 4. For example,
Theorem~\ref{theorem:comulative_small} gives the exact values of
$e(3,9,n)$ for all $n \le 26$, of $e(3,10,n)$ for $n \le 28$,
and of $e(3,13,n)$ for all $n \le 39$.

\smallskip
The inequality $e(3,k+1,n) \ge (40n-91k)/6$, which is better than
$e(3,k+1,n) \ge 6n-13k$ for larger parameters,
and a number of other improvements and characterizations
of graphs realizing specific number of edges, was credited
in 2001 by Lesser \cite{Les} to an unpublished manuscript
by Backelin \cite{Back}. As of 2012, the manuscript by Backelin
already exceeds 500 pages and it contains numerous additional
related results \cite{Back, BackP}, but it still needs more
work before it can be published. Therefore, in the remainder
of this paper we will not rely on the results included therein,
however in several places we will cite the bounds obtained there
for reference. In summary, the behavior of $e(3,k+1,n)$ is clear
for $n \le 13k/4-1$, it seems regular but very difficult to deal
with for $n$ slightly larger than $13k/4$, and
it becomes hopelessly hard for even larger $n$.
In this work we apply computational techniques
to establish lower bounds for $e(3,k,n)$ for larger $n$, for
$k \le 15$. Immediately,
our results imply better upper bounds on $R(3,k)$ in several cases,
but we hope that they also may contribute to further
progress in understanding the general behavior of $e(3,k,n)$.

\medskip
Full enumeration of the sets $\mathcal{R}(3,\le 6)$ was established
in \cite{RK1,MZ}. The knowledge of the exact values of $e(3,7,n)$
was completed in \cite{RK1}, those of $e(3,8,\le 26)$ in \cite{RK2},
and the last missing value for $\alpha(G)<8$, namely $e(3,8,27)=85$,
was obtained in \cite{BGSP}.
The thesis by Lesser \cite{Les} contains many lower bounds
on $e(3,k,n)$ better than those in \cite{RK2}. We match or improve
them in all cases for $k \le 10$. For
$k \ge 11$ and $n$ slightly exceeding $13k/4-1$, the bounds
by Lesser (in part credited also to \cite{Back}) are better
than ours in several cases, however we obtain significantly
better ones for larger $n$.

\medskip
The general method we use is first to compute, if feasible,
the exact value of $e(3,k,n)$ for concrete $k$ and $n$,
or to derive a lower bound using
a combination of (2), (3) and (4), and computations. Better lower
bounds on $e(3,k-1,m)$ for $m=n-d-1$ and various $d$, lead in
general to better lower bounds on $e(3,k,n)$. If we manage
to show that $e(3,k,n)=\infty$, i.e. no $(3,k;n)$-graph exists,
then we obtain an upper bound $R(3,k)\le n$. An additional
specialized algorithm was needed to establish $R(3,10)\le 42$.

\bigskip
Section~\ref{section:algorithms} describes extension algorithms
which we used to exhaustively construct all $(3,k;n,e)$-graphs
for a number of cases of $(n,e)$, for $k \le 10$.
These results are described in
detail in the sequel. This leads to many new lower bounds
on $e(3,k,n)$ and full enumerations of $(3,k;n)$-graphs with
the number of edges equal to or little larger than $e(3,k,n)$,
which are presented in Section~\ref{section:computations_10}
(and Appendix 1).
These results are then used in Section~\ref{section:better_bounds}
to prove that there exists a unique critical
35-vertex graph for the Ramsey number $R(3,9)$.
It is known that \cite{Ex} $40 \le R(3,10)\le 43$ \cite{RK2}.
We establish that $R(3,10)=43$ if and only if $e(3,10,42)=189$,
or equivalently, that if $R(3,10)=43$ then every critical graph
in this case is regular of degree 9. Then, in Section 6,
using computations we prove that the latter do not exist,
and thus obtain $R(3,10)\le 42$.
Finally, in Section~\ref{section:larger_10}, we describe
the second stage of
our computations, which imply many new lower bounds on
$e(3,\ge 11,n)$. This stage uses only degree sequence analysis
of potential $(3,k;n,e)$-graphs, which have to
satisfy (4). This in turn leads to the
new upper bounds on the classical two-color
Ramsey numbers marked in bold in Table~\ref{table:new_bounds},
which presents the values and best bounds on the Ramsey
numbers $R(3,k)$ for $k \le 16$. All the improvements
in this work are better by one over the results listed in
the latest 2011 revision \#13 of the survey \cite{SRN}.
The bound $R(3,16) \le 98$ was also obtained by Backelin
in 2004, though it was not published \cite{Back,BackP}.
The lower bound $R(3,11) \ge 47$ was recently obtained
by Exoo \cite{ExP}. The references for all other bounds
and values, and the previous upper bounds, are
listed in \cite{SRN}.

\bigskip
\begin{table}[H]
\begin{center}
\begin{tabular}{|c|c||c|c|}
\hline
$k$&$R(3,k)$&$k$&$R(3,k)$\cr
\hline
3&\ \ 6&10&40--{\bf 42}\cr
4&\ \ 9&11&47--{\bf 50}\cr
5&14&12&52--59\cr
6&18&13&59--{\bf 68}\cr
7&23&14&66--{\bf 77}\cr
8&28&15&73--{\bf 87}\cr
9&36&16&79--{\bf 98}\cr
\hline
\end{tabular}

\caption{Ramsey numbers $R(3,k)$, for $k \le 16$.}
\label{table:new_bounds}
\end{center}

\end{table}

\section{Algorithms}
\label{section:algorithms}

\subsection*{Maximum Triangle-Free Method}

One method to determine $e(3,k,n)$ is by first generating all \textit{maximal}
triangle-free $(3,k;n)$-graphs. A maximal triangle-free graph (in short,
an \textit{mtf graph}) is a triangle-free graph such that the insertion of
any new edge forms a triangle. It is easy to see that there exists a
$(3,k;n)$-graph if and only if there is an mtf $(3,k;n)$-graph.
In \cite{BGSP}, an algorithm is described that can generate all mtf
$(3,k;n)$-graphs efficiently. Using this algorithm, it is much easier
to generate all mtf $(3,k;n)$-graphs instead of all $(3,k;n)$-graphs,
because the number of the former is in most cases much smaller.
For example, there are 477142 $(3,8;27)$-graphs, but only
21798 mtf graphs with the same parameters. By recursively removing
edges in all possible ways from these mtf $(3,k;n)$-graphs and testing
if the resulting graphs $G$ still satisfy $\alpha(G)<k$,
the complete set $\mathcal{R}(3,k;n)$ can be obtained.

We applied this method to generate the sets
$\mathcal{R}(3,7;21), \mathcal{R}(3,7;22), \mathcal{R}(3,8;26,\le 77)$
and $\mathcal{R}(3,8;27)$ (see Appendix 1 for detailed results).
All $(3,7;22)$- and $(3,7;n,e(3,k,n))$-graphs were already
known \cite{RK1}, other enumerations are new. This mtf method
is infeasible for generating $(3,\ge 9;n)$-graphs for $n$ which
were needed in this work. Nevertheless, we used it for verifying
the correctness of our other enumerations, and the results agreed in
all cases in which more than one method was used (see Appendix 2).

\subsection*{Minimum Degree Extension Method}

In their 1992 paper establishing $R(3,8)=28$,
McKay and Zhang \cite{MZ}
proved that the set $\mathcal{R}(3,8;28)$ is empty
by generating several sets $\mathcal{R}(3,k;n,e)$ with additional
restrictions on the minimum degree $\delta(G)$.
Suppose that one wants to
generate all $(3,k;n,e)$-graphs. If $G$ is such a graph and one
considers its minimum degree vertex $v$, then we can reconstruct
$G$ given all possible graphs $G_v$. McKay and Zhang described
such dependencies, designed an algorithm to reconstruct $G$,
and completed the proof of $R(3,8)=28$ using this algorithm.

We implemented and used this method by McKay and Zhang \cite{MZ},
and in all cases where more than one algorithm was used it agreed
with the other results. However, using this method it was not feasible
to generate most classes of graphs with higher parameters needed
for our project. For example,
we could not generate all $(3,9;28,\le 69)$-graphs with this method,
as the graphs with $\delta(G)=4$ are obtained from
$(3,8;23,\le 53)$-graphs, but there are already 10691100
$(3,8;23,\le 52)$-graphs (Table~\ref{table:graph_counts_K8}
in Appendix 1).

\subsection*{Neighborhood Gluing Extension Method}

Our general extension algorithm for an input $(3,k;m)$-graph $H$
produces all $(3,k+1;n,e)$-graphs $G$, often with some specific
restrictions on $n$ and $e$, such that for some vertex $v \in V(G)$
graph $H$ is isomorphic to $G_v$.
We used the following strategy to determine if the parameters
of input graphs to our extender were such that the output was
guaranteed to contain all $(3,k+1;n,\le e)$-graphs.

Let $m_i=n-i-1$, where $i$ ranges over possible degrees
in any graph $G$ we look for, $\delta(G) \le i \le \Delta(G)$.
In the broadest case we have $\delta(G)=\max\{n-R(3,k),0\}$
and $\Delta(G)=k$, but we also identified a number
of special cases where this range was more restricted.
Let $t_i$ be an integer such that we have extended
all $(3,k;m_i,< e(3,k,m_i)+t_i)$-graphs as potential
$G_v$'s of $G$. Now, if we use $e(3,k,m_i)+t_i$ instead
of $e(3,k,m_i)$ in (4) for all relevant values of $i$,
and (4) has no solutions for $(3,k+1;n,\le e)$-graphs,
then we can conclude that all such graphs were already
generated. We illustrate this process by an example.

\bigskip
\noindent
{\bf Example.}
Table~\ref{table:example_glueing} lists specific parameters
of the general process when used to obtain all
$(3,8;25, \le 65)$-graphs.
Every vertex $v$ in any $(3,8;25, \le 65)$-graph has degree
$i$, for some $2 \le i \le 7$. The corresponding graph
$G_v$ is of type $(3,7;m_i,e(G_v))$. The values
of $e(3,7,m)$ are included in Table~\ref{table:exact_e_values_all}
of Section~\ref{section:computations_10},
and let $t_i$'s be as in Table 2.
If we use the values $e(3,7,m_i)+t_i$ instead of $e(3,7,m_i)$
in (4), then there are no solutions for degree sequences
of $(3,8;25, \le 65)$-graphs. Thus, if we run the
extender for all possible graphs $G_v$ with the number
of edges listed in the last column of Table~\ref{table:example_glueing},
then we will obtain all $(3,8;25,e)$-graphs for $e \le 65$.

The set of increments $t_i$ accomplishing this goal is not
unique, there are others which work. We just tried to minimize
the amount of required computations in a greedy way.
Note that the largest increments $t_i$ to $e(3,7,m_i)$ occur
for $i$'s which are close to the average degree of $G$.

\medskip
\begin{table}[H]
\begin{center}
\begin{tabular}{|c|c|c|l|c|}
\hline
$i=\deg_G(v)$&$m_i=|V(G_v)|$&$e(3,7,m_i)$&$t_i$&$e(G_v)=e-Z(v)$\cr
\hline
2&22&60&1&60\cr
3&21&51&1&51\cr
4&20&44&2&44, 45\cr
5&19&37&3&37, 38, 39\cr
6&18&30&2&30, 31\cr
7&17&25&1&25\cr
\hline
\end{tabular}


\caption{Obtaining all $(3,8;25, \le 65)$-graphs.}
\label{table:example_glueing}
\end{center}
\end{table}

\subsection*{Implementation}

In this section we present some details about the extension algorithms
implementations for the minimum degree and neighborhood gluing method.
Implementation of the algorithm  to generate maximal triangle-free
Ramsey graphs is described in~\cite{BGSP}.

Given a $(3,k;n,f)$-graph $G'$ as input and an expansion degree $d$,
a desired maximum number of edges $e$, and the minimum degree $d_m$ as
parameters, our program constructs all $(3,k+1;n+d+1,\le e)$-graphs
$G$ with $\delta(G) \ge d_m$ for which there is a vertex $v \in V(G)$
such that $\deg(v) = d$ and $G_v \cong G'$.
More specifically, the program adds to $G'$ a vertex $v$ with neighbors
$u_1,...,u_d$ and connects them to independent sets of $G'$ in all
possible ways, so that the resulting graph is a
$(3,k+1;n+d+1,\le e)$-graph with $\delta(G) \ge d_m$.
Note that the neighbors of $v$ have to be connected to
independent sets of $G'$, otherwise the expanded graph
would contain triangles, and, clearly, $\Delta(G)\le k$.

The extension program first determines all independent sets of
$G'$ of orders $t$ that are possible, namely $d_m - 1 \le t \le k-1$.
The program then recursively assigns the $d$ neighbors
of $v$ to the eligible independent sets of $G'$,
adds the edges joining $u_i$'s to their associated independent sets,
and tests if the resulting $G$ is a valid $(3,k+1;n+d+1,\le e)$-graph.
If it is, then we output it. This general process is greatly accelerated
by the techniques described in the following.

We bound the recursion if a given partial assignment cannot lead to
any $(3,k+1;n+d+1,\le e)$-graphs. Suppose that $i$ independent sets
$S_1,\dots ,S_i$ have already been assigned.
If $V(G') \setminus (S_1 \cup ... \cup S_i)$ induces an independent
set $I$ of order $k + 1 - i$, then
this assignment cannot lead to any output since
$I \cup \{u_1,\dots ,u_i\}$ would form an independent set
of order $k+1$ in $G$.
We could test this property for all subsets of $S_i$'s,
but we found it to be most efficient to do it only for all pairs.
Namely, if
$S_1,\dots ,S_i$ is already assigned and we consider the next
independent set $S$, we test if for all $j$, $1 \le j \le i$,
$V(G') \setminus (S_j \cup S)$ does not induce
any independent set of order $k-1$.
The list of independent sets which can still
be assigned is dynamically updated.

For the efficiency of the algorithm it is vital that testing
for independence in $V(G') \setminus (S_1 \cup ... \cup S_i)$ 
is fast, and hence we precompute the independence numbers of
all induced subgraphs of $G'$.
This precomputation also needs to be done very efficiently. 
We represent a set of vertices $S \subset V(G')$ by a bitvector.
The array
\verb|indep_number[S]| of $2^n$ elements stores the independence
number of the graph induced by $S$ in $G'$.
It is very important that \verb|indep_number[]| fits
into the memory. On the computers on
which we performed the expansions this was still feasible up to $n=31$.
We investigated various approaches to precompute \verb|indep_number[S]|,
and Algorithm 1 below was by far be the most efficient one.
If the superset $S'$ of $S$ already has
\verb|indep_number[S']| $\ge j$, then we can break the recursion of making
the supersets. Usually one can break very quickly.
For small extension
degrees $d \le 3$, it is more efficient not to precompute these
independence numbers, but instead to compute them
as needed.

\medskip
\begin{algorithm}[H]
\caption{Precomputing independence number}
  \begin{algorithmic}
  \label{algo:compute_indep_numbers}
	\FOR{$i = 0$ {\bf upto} $2^n-1$}
		\STATE set \verb|indep_number[i]| $=0$
	\ENDFOR
  	\FOR{$j= k - 1$ {\bf downto} $k + 1 - d$}
	\FOR{all independent sets $S$ of order $j$ in $G'$}
		\STATE Recursively make all supersets $S'$ of $S$, and\\
		{\bf if} \verb|indep_number[S']| $=0$ {\bf then} set \verb|indep_number[S']| $=j$\\
		{\bf else} break making supersets of $S$
	\ENDFOR
	\ENDFOR
  \end{algorithmic}
\end{algorithm}

\medskip
If a neighbor $u_i$ of $v$ has been assigned to an independent set $S$,
we also update the degrees of the vertices in $G'$. If $u_i$ is being connected
to $S$, the degree of every vertex of $S$ increases by one. If the degree
of a vertex $w$ of $G'$ becomes $k$, then other neighbors
of $v$ cannot be assigned to independent sets which contain $w$.
We call such vertices which are no longer eligible \textit{forbidden vertices},
and all of them are stored in a dynamically updated bitvector.
We also dynamically update the list of independent sets which can still
be assigned to $u_i$'s. Independent sets which contain forbidden vertices
are removed from the list of eligible independent sets.
We perform bitvector operations whenever suitable.
If no eligible independent sets are left, we can bound the recursion.
Note that we cannot break the recursion when the number of eligible
independent sets is smaller than the number of neighbors of $v$ that still
have to be considered, since they can be assigned to the same independent set.
If $i$ neighbors of $v$ are already assigned and the forbidden vertices
form an independent of set order $k + 1 - (d - i)$, then the recursion
can also be bounded, though this criterion in general is weak.

We assign the neighbors $u_i$ of $v$ to independent
sets in ascending order, i.e. if $u_i$ is assigned to $S_i$,
then $|S_i| \le |S_{i+1}|$ for all $1 \le i < d$.
Doing this rather than in descending order allows us to eliminate
many candidate independent sets early in the recursion.
If $|S_i|$ is small, then it is very likely that
$V(G') \setminus S_i$ induces a large independent set.
Hence, it is also very likely that $S_i$ cannot be assigned
to a new $u_i$ or that assigning $S_i$ eliminates many eligible
independent sets.

Assigning sets in ascending order also gives us an easy lower bound
for the number of edges in any potential output graph which can be
obtained from the current graph and assignment. If the sets
$S_1,\dots ,S_i$ have already been
assigned to neighbors of $v$ and the current minimal order of
eligible independent sets is $t$, then any expanded graph will have at least
$f = e(G') + d + |S_1| + ... + |S_i| + t(d - i)$ edges. If $f > e$,
then we can bound the recursion as well.

The pseudocode of the recursive extension is listed below
as Algorithm 2. It is assumed that \verb|indep_number[]|
(see Algorithm~\ref{algo:compute_indep_numbers}) and
the list of eligible independent sets are already computed.
The parameters for \verb|Construct()| are the order of
the sets which are currently being assigned and the number
of neighbors of $v$ which were already assigned to independent sets.
The recursion is bounded if any of the bounding criteria described
above can be applied.

\medskip
\begin{algorithm}[H]
\caption{{\tt Construct}(current\_order, num\_assigned)}
  \begin{algorithmic}
	\IF{num\_assigned $=d$}
  		\STATE expand graph $G'$ to $G$
  		\IF{$G$ is a $(3,k+1;n+d+1,\le e)$-graph}  
  			\STATE output $G$
  		\ENDIF
  	\ELSE
  		\FOR{every eligible set $S$ of order current\_order}
			\STATE assign $S$ to $u_{num\_assigned+1}$
			\STATE update the set of eligible independent sets
			\STATE \verb|Construct|(current\_order, num\_assigned + 1)
		\ENDFOR
		\IF{current\_order $< k - 1$}  
			\STATE \verb|Construct|(current\_order + 1, num\_assigned)
		\ENDIF
  	\ENDIF
  \end{algorithmic}
\end{algorithm}

\smallskip
Our extension program does not perform any isomorphism rejection.
We canonically label the output graphs with
\textit{nauty}~\cite{McKay, nauty}
and remove the isomorphic copies.
This is not a bottleneck as there are usually only a few
$(3,k+1;n+d+1,\le e)$-graphs which are constructed by our program.
The results obtained by our extension algorithms are described in
Sections~\ref{section:computations_10} and 6.
In the appendices we describe
how the correctness of our implementation was tested.

\subsection*{Degree Sequence Feasibility}

Suppose we know the values or lower bounds on $e(3,k,m)$ for some fixed $k$
and we wish to know all feasible degree sequences of $(3,k+1;n,e)$-graphs.
We construct the system of integer constraints consisting of
$n=\sum_{i=0}^{k}{n_i}$, $2e=\sum_{i=0}^{k}{i n_i}$, and (4).
If it has no solutions then we conclude that $e(3,k+1,n) > e$.
Otherwise, we obtain solutions for $n_i$'s which include all
desired degree sequences.
This algorithm is similar in functionality to the package FRANK developed
by Lesser \cite{Les}.

\section{Progress on Computing Small $e(3,k,n)$}
\label{section:computations_10}

\begin{table}[H]
\begin{center}
\begin{tabular}{|c|rrrrrrrrrrrrrr|}
\hline
vertices&\multicolumn{14}{c}{$k$}\vline\cr
$n$&3&4&5&6&7&8&9&10&11&12&13&14&15&16\cr
\hline
3&1&&&&&&&&&&&&&\cr
4&2&1&&&&&&&&&&&&\cr
5&5&2&1&&&&&&&&&&&\cr
6&$\infty$&3&2&1&&&&&&&&&&\cr
7&&6&3&2&1&&&&&&&&&\cr
8&&10&4&3&2&1&&&&&&&&\cr
9&&$\infty$&7&4&3&2&1&&&&&&&\cr
10&&&10&5&4&3&2&1&&&&&&\cr
11&&&15&8&5&4&3&2&1&&&&&\cr
12&&&20&11&6&5&4&3&2&1&&&&\cr
13&&&26&15&9&6&5&4&3&2&1&&&\cr
14&&&$\infty$&20&12&7&6&5&4&3&2&1&&\cr
15&&&&25&15&10&7&6&5&4&3&2&1&\cr
16&&&&{\bf 32}&20&13&8&7&6&5&4&3&2&1\cr
17&&&&{\bf 40}&25&16&11&8&7&6&5&4&3&2\cr
18&&&&$\infty$&30&20&14&9&8&7&6&5&4&3\cr
19&&&&&{\bf 37}&25&17&12&9&8&7&6&5&4\cr
20&&&&&{\bf 44}&30&20&15&10&9&8&7&6&5\cr
21&&&&&{\bf 51}&35&25&18&13&10&9&8&7&6\cr
22&&&&&{\bf 60}&{\bf 42}&30&21&16&11&10&9&8&7\cr
23&&&&&$\infty$&{\bf 49}&35&25&19&14&11&10&9&8\cr
24&&&&&&{\bf 56}&40&30&22&17&12&11&10&9\cr
25&&&&&&{\bf 65}&46&35&25&20&15&12&11&10\cr
26&&&&&&{\bf 73}&52&40&30&23&18&13&12&11\cr
27&&&&&&{\bf 85}&{\bf 61}&45&35&26&21&16&13&12\cr
28&&&&&&$\infty$&{\bf 68}&51&40&30&24&19&14&13\cr
29&&&&&&&{\bf 77}&{\bf 58}&45&35&27&22&17&14\cr
30&&&&&&&{\bf 86}&{\bf 66}&50&40&30&25&20&15\cr
31&&&&&&&{\bf 95}&{\bf 73}&56&45&35&28&23&18\cr
\hline
\end{tabular}

\caption{Exact values of $e(3,k,n)$, for $3 \le k \le 16$, $3 \le n \le 31$.}
\label{table:exact_e_values_all}
\end{center}
\end{table}

Table~\ref{table:exact_e_values_all} presents the exact values of $e(3,k,n)$ for small cases,
where clear regularities are well described by Theorem~\ref{theorem:comulative_small}. Empty
entries in the upper-right triangle of the table are 0's,
while those in the lower-left
triangle are equal to $\infty$. The columns correspond to fixed
values of $k$. Almost all entries are given by Theorem~\ref{theorem:comulative_small}.
We list them for a better perspective and completeness.
The entries beyond the range of Theorem~\ref{theorem:comulative_small}
are marked in bold, and they were obtained as follows:
$e(3,6,16)$ and $e(3,6,17)$ in \cite{GY},
all cases for $k=7$ in \cite{GY, GR, RK1},
all cases for $k=8$ and $22 \le n \le 26$ in \cite{RK2},
$e(3,8,27)=85$ was computed in \cite{BGSP},
and those for $k \ge 9$ are obtained here.
The smallest $n$ for which we found an open case is 32,
namely that of $e(3,11,32)$. However, Backelin had claimed
the values $e(3,11,32)=63$, $e(3,11,33)=70$ and $e(3,11,34)=77$ \cite{Back,BackP}.
Tables~\ref{table:exact_e_values_9} and~\ref{table:bounds_e_10}
below and 7--11 in Section~\ref{section:larger_10} present
the details of what we found about these harder parts of
each column $k$, for $9 \le k \le 16$.

The exact counts of $(3,k;n,e)$-graphs for $k=7,8,9,10$ which were obtained
by the algorithms described in Section~\ref{section:algorithms} are listed in
Tables 12, 13, 14, 15, respectively, in Appendix 1.
All $(3,\le 9;n,\le e(3,k,n)+1)$-graphs which were constructed by our programs can
be obtained from the \textit{House of Graphs}~\cite{HOG} by searching for
the keywords ``minimal ramsey graph'' or from \cite{hog-static}.

\medskip
\bigskip
\noindent
{\bf Exact values of $e(3,9,n)$}

\medskip
\noindent
The values of $e(3,9,\le 26)$ are determined by
Theorem~\ref{theorem:comulative_small}.
The values of $e(3,9,n)$ for $27 \le n \le 34$ were obtained
by computations, mostly by the gluing extender algorithm described
in Section~\ref{section:algorithms}, and they are presented
in Table~\ref{table:exact_e_values_9}.
All of these values improve over previously reported lower
bounds \cite{RK2,Les}. The equality $e(3,9,35)=140$ will
be established by Theorem~\ref{theorem:r_3_9_35}
in Section~\ref{section:better_bounds}.

\begin{table}[H]
\begin{center}
\begin{tabular}{c|c|l}
\hline
$n$&$e(3,9,n)$&comments\cr
\hline
27&\ \ 61&\cr
28&\ \ 68&the same as in \cite{Back}\cr
29&\ \ 77&\cr
30&\ \ 86&\cr
31&\ \ 95&\cr
32&104&not enough for $R(3,10) \le 42$\cr
33&118&just enough for Theorem~\ref{theorem:r_3_10_42}\cr
34&129&122 required for $R(3,10) \le 43$\cr
35&140&Theorem~\ref{theorem:r_3_9_35}\cr
36&$\infty$&hence $R(3,9) \le 36$, old bound\cr
\hline
\end{tabular}

\caption{Exact values of $e(3,9,n)$, for $n\ge 27$}
\label{table:exact_e_values_9}
\end{center}

\end{table}

\noindent
{\bf Values and lower bounds on $e(3,10,n)$}

\medskip
\noindent
The values of $e(3,10,\le 28)$ are determined by Theorem~\ref{theorem:comulative_small}.
The values for $29 \le n \le 34$ were obtained
by the gluing extender algorithm described in Section~\ref{section:algorithms}.
The lower bounds on $e(3,10,\ge 35)$
are included in the second column of Table~\ref{table:bounds_e_10}.
They are based on solving
integer constraints (3) and (4), using the exact values
of $e(3,9,n)$ listed in Table 4,
and results from the gluing extender algorithm used similarly as
in the example of Section 3.
Our bounds on $e(3,10,n)$ improve over previously
reported lower bounds \cite{RK2,Les} for all $n \ge 30$.

\medskip
By Theorem~\ref{theorem:r_3_10_42}
(see Section~\ref{section:better_bounds}) we know that
any $(3,10;42)$-graph must be 9-regular
with 189 edges, and thus all its graphs $G_v$ are necessarily of
the type $(3,9;32,108)$. There exists a very large
number of the latter graphs. Their generation, extensions
to possible $(3,10;42,189)$-graphs, and implied nonexistence
of any $(3,10;42)$-graphs will be described in Section 6.

\bigskip
\begin{table}[H]
\begin{center}
\begin{tabular}{c|c|l}
\hline
$n$&$e(3,10,n)\ge$&comments\cr
\hline
29&\ \ 58&exact, the same as in \cite{Les}\cr
30&\ \ 66&exact\cr
31&\ \ 73&exact\cr
32&\ \ 81&exact\cr
33&\ \ 90&exact\cr
34&\ \ 99&exact, $(3,10;34,99)$-graph constructed by Backelin \cite{BackP}\cr
35&107\cr
36&117\cr
37&128\cr
38&139&146 required for $R(3,11) \le 49$\cr
39&151&as required for $R(3,11) \le 50$, Theorem~\ref{theorem:new_upper_bounds}\cr
40&161&\cr
41&172&184 maximum\cr
42&$\infty$&hence $R(3,10) \le 42$, new bound, Theorem 6\cr
43&$\infty$&hence $R(3,10) \le 43$, old bound\cr
\hline
\end{tabular}

\caption{Values and lower bounds on $e(3,10,n)$, for $n\ge 29$.}
\label{table:bounds_e_10}
\end{center}
\end{table}

All lower bounds in Tables 4 and 5 required computations of
our gluing extender algorithm. We did not perform any such
computations in an attempt to improve the lower bounds on
$e(3,\ge 11,n)$. All results presented in Section 7 for
$k \ge 11$ depend only on the degree sequence analysis
and the results for $k \le 10$.

\medskip

\section{Better Lower Bounds for $e(3,9,35)$ and $e(3,10,42)$}
\label{section:better_bounds}

Sometimes we can improve on the lower bounds on $e(3,k,n)$
implied by (3) and (4) by a more detailed analysis of
feasible degree sequences. Such improvements typically can
be done in cases for which (4) gives a small number of possible
degree sequences, none of which is of a regular graph,
furthermore with only one heavily dominating
degree. We have such a situation in
the proofs of the two following theorems.  

\begin{thm}
\label{theorem:r_3_9_35}
There exists a unique $(3,9;35)$-graph, and $e(3,9,35)=140$.
\end{thm}

\begin{proof}
Any $(3,9;35)$-graph $G$ has $\Delta(G)\le 8$,
hence we have $e(G)\le 140$.
Suppose $G \in \mathcal{R}(3,9;35,140-s)$ for some $s\ge 0$.
Since $R(3,8)=28$, the degrees of vertices in $G$
are 7 or 8, and let there be $n_7$ and $n_8$ of them,
respectively. We have $n_7+n_8=35$, $n_7=2s$. In this
case there are five solutions to (4) with $0 \le s \le 4$.
In particular, this shows that $e(3,9,35)\ge 136$.
If $n_7>0$ (equivalently $s>0$), then
consider graph $H$ induced in $G$ by $n_7$ vertices of degree 7.
Observe that $\delta(H)\le s$, since $H$ is triangle-free on $2s$ vertices.
Let $v$ be a vertex in $V(G)$ of degree 7 connected to at most $s$
other vertices of degree 7. Thus we have $Z_G(v) \ge 7s+8(7-s) = 56 -s$,
and $e(G_v) \le (140-s)-(56-s)=84$. However $G_v$ is a $(3,8;27)$-graph
which contradicts the fact that $e(3,8,27)=85$.

The computations extending all $(3,8;26,76)$-graphs,
using the neighborhood gluing extension method
described in Section~\ref{section:algorithms}, established
that there exists a unique (up to isomorphism) 8-regular
$(3,9;35)$-graph.
We note that it is a cyclic graph on 35 vertices with circular
distances \{1,7,11,16\}, found by Kalbfleisch~\cite{Kalb} in 1966.
Clearly, any  $(3,9;35,140)$-graph must be 8-regular,
and thus the theorem follows.
\end{proof}

\bigskip
\begin{thm}
\label{theorem:r_3_10_42}
$R(3,10)=43$\ \ \  if and only if\ \ \  $e(3,10,42)=189$.
\end{thm}

\begin{proof}
It is known that $R(3,10) \le 43$ \cite{RK2},
i.e. there are no $(3,10;43)$-graphs.
We will prove the theorem by showing that any $(3,10;42)$-graph
must be regular of degree 9. The essence of the reasoning is very
similar to that for $e(3,9,35)=140$ in the previous theorem,
except that this time it is little more complicated.

Suppose $G \in \mathcal{R}(3,10;42,189-s)$ for some $s\ge 0$.
The computations described in Section~\ref{section:algorithms} established
that $G$ cannot have the unique $(3,9;35)$-graph as one of its $G_v$'s.
Hence, $7 \le \deg_G(v) \le 9$ for all vertices $v \in V(G)$.
The solutions $n_i$ to (4) which contain all possible degree sequences
for $G$  with this restriction are presented in Table~\ref{table:solutions_r310}.

\begin{table}[H]
\begin{center}
\begin{tabular}{ccc|ccc}
\hline
$n_7$&$n_8$&$n_9$&$e(G)$&$\gamma(G)$&$s$\cr
\hline
\hline
0&8&34&185&24&4\cr
1&6&35&185&25&4\cr
2&4&36&185&26&4\cr
3&2&37&185&27&4\cr
4&0&38&185&28&4\cr
\hline
0&6&36&186&60&3\cr
1&4&37&186&61&3\cr
2&2&38&186&62&3\cr
3&0&39&186&63&3\cr
\hline
0&4&38&187&96&2\cr
1&2&39&187&97&2\cr
2&0&40&187&98&2\cr
\hline
0&2&40&188&132&1\cr
1&0&41&188&133&1\cr
\hline
0&0&42&189&168&0\cr
\hline
\end{tabular}

\caption{Solutions to (4) for $(3,10;42,189-s)$-graphs.}
\label{table:solutions_r310}
\end{center}
\end{table}

\smallskip
Note that for all $0 \le s \le 4$ we have $0 \le n_7 \le s$,
$n_8+2n_7=2s$, $n_9=42-n_8-n_7$, and $e(G)=189-s$.
Since $e(3,9,34)=129$, using (2) we see that
$Z(v) \le 60-s$ for every vertex $v$ of degree 7.
Similarly, since $e(3,9,33)=118$,
$Z(v) \le 71-s$ for every vertex $v$ of degree 8.
If $s=0$, then we are done, otherwise
consider graph $H$ induced in $G$ by $2s-n_7$ vertices of degree 7 or 8.
Observe that $\delta(H)\le s -n_7/2$, since $H$ is triangle-free.

\medskip
\noindent
{\bf Case 1: $n_7=0$.}
Let $v$ be a vertex in $V(G)$ of degree 8 connected to at most $s$
other vertices of degree 8. This gives
$Z_G(v) \ge 8s+9(8-s) = 72-s$, which is a contradiction.

\medskip
\noindent
{\bf Case 2: $n_8=0$.}
Let $v$ be a vertex in $V(G)$ of degree 7 connected to at most $s/2$
other vertices of degree 7 (in this case $|V(H)|=s$). This gives
$Z_G(v) \ge 7s/2+9(7-s/2) = 63-s$, which is a contradiction.

\medskip
\noindent
{\bf Case 3: $n_7=1$.}
If $v$ is the only vertex of degree 7, then $n_8=2s-2$ and
we easily have $Z_G(v) \ge 8n_8 + 9(7-n_8) = 65-2s>60-s$, which
again is a contradiction.

\medskip
\noindent
{\bf Case 4: $n_7=2$.}
Both vertices of degree 7 must have
$Z_G(v) \ge 7 + 8n_8 + 9(7-n_8-1)=61-(2s-2n_7)=65-2s$,
which is a contradiction.

\medskip
\noindent
{\bf Case 5: $n_7>2$.} The only remaining degree sequence not
covered by previous cases is
$n_7=3$ and $n_8=2$, for $s=4$ and $e=185$.
There is a vertex $v$ of degree 7 connected to at most one
other of degree 7, and thus
$Z_G(v) \ge 7+2\cdot 8+4\cdot 9 > 60-s$, a contradiction.
\end{proof}

\bigskip

\section{$R(3,10) \le 42$}

\medskip
Theorem 4 implies that any $(3,10;42)$-graph $G$ must be regular of
degree 9 with 189 edges. Removing any vertex
$v$ with its neighborhood from $G$ yields a $(3,9;32,108)$-graph $G_v$.
Hence, our first task is to obtain all $(3,9;32,108)$-graphs.

We used the neighborhood extension method to
generate $(3,9;32,108)$-graphs $H$ with a vertex
$v$ for which $H_v$ is one of the following types:
$(3,8;27)$, $(3,8;26,\le 77)$, $(3,8;25,\le 68)$,
$(3,8;24,\le 59)$ or $(3,8;23,49)$.
These extensions yielded the set of 2104151
$(3,9;32,108)$-graphs $\mathcal{X}$.
Using notation of the example in Section 3, now with
$4 \le i \le 8$, $m_i=31-i$, and $t_i=10,5,4,4,1$,
respectively, the only remaining degree sequence
passing (4) for a $(3,9;32,108)$-graph
is $n_6=8, n_7=24$.

Potentially, the complete set of $(3,9;32,108)$-graphs
could be obtained by performing additional extensions of
degree 6 to $(3,8;25,69)$-graphs or extensions of
degree 7 to $(3,8;24,60)$-graphs.
However, there are already
12581543 $(3,8;25,\le 68)$-graphs and 3421512
$(3,8;24,\le 59)$-graphs (see Table 13 in Appendix 1),
and there are many more with one additional edge.
Hence, further refinement of the construction method
of the $(3,9;32,108)$-graphs not in $\mathcal{X}$
was needed. It is described in the following 
Lemma~\ref{theorem:restricted_n6_8__n7_24},
which permitted a fast computation and
the completion of the task.

\medskip
\begin{lem}
\label{theorem:restricted_n6_8__n7_24}
Every $(3,9;32,108)$-graph $H \not\in \mathcal{X}$ has
$n_6=8, n_7=24$, and furthermore in such $H$
every vertex of degree $6$ has exactly $3$ neighbors
of degree $7$ and every vertex of degree $7$ has exactly $1$
neighbor of degree $6$.
\end{lem}

\begin{proof}
As stated after the definition of $\mathcal{X}$ above,
(4) implies the specified degree sequence of $H \not\in \mathcal{X}$.
Suppose that $H$ has a vertex $v$ of degree 6 with at least 4
neighbors of degree 7. One can easily see that
$Z_H(v)\ge 40$ and thus $e(H_v)\le 68$. All such graphs, however,
were included in the set of inputs producing $\mathcal{X}$,
so we have a contradiction. Similarly,
suppose that $H$ has a vertex $v$ of degree 7 with no
neighbors of degree 6. Then  $Z_H(v)=49$ and $e(H_v)=59$,
but all such graphs were used as inputs producing $\mathcal{X}$,
hence again we have a contradiction. Now,
by the pigeonhole principle, there are exactly 24 edges
connecting vertices of distinct degrees, and we can easily
conclude that every vertex of degree 6
must have exactly 3 neighbors of degree 7 and every
vertex of degree 7 exactly 1 neighbor of degree 6.
\end{proof}

We adapted the extension algorithm from Section 3 to generate
this very restricted set of $(3,9;32,108)$-graphs by performing
extensions of all 64233886 $(3,8;24,60)$-graphs (Table 13
in Appendix 1). The result is that there are
no $(3,9;32,108)$-graphs not in $\mathcal{X}$.

\bigskip
\begin{thm}
\label{theorem:r_3_10_42b}
$R(3,10)\le 42$.
\end{thm}

\begin{proof}
For contradiction, suppose that $G$ is a $(3,10;42)$-graph.
By Theorem 4 it must be a 9-regular
$(3,10;42,189)$-graph whose all $G_v$'s are $(3,9;32,108)$-graphs.
By Lemma 5 and the computations described above there are exactly
2104151 such graphs. A specialized extension algorithm (a modification
of the gluing extender) was run
for all of them in an attempt to obtain a 9-regular $(3,10;42,189)$-graph.
The neighbors of $v$ have to be connected to independent sets of
order 8 in $G_v$. For every pair of (possibly equal)
independent sets $\{S_i, S_j\}$ of order 8, we
test if they can be assigned to two neighbors of $v$ by checking
if $V(G_v) \setminus (S_i \cup S_j)$ induces an independent
set of order 8 in $G_v$, and if so we can bound the recursion.
We used for this task a precomputed table storing the results
of such tests for all pairs of independent sets of order 8.
The concept of eligible candidates (Section 3) was also used,
and the condition $\Delta(G)=9$ turned out to be particularly
strong in pruning the recursion.
No 9-regular $(3,10;42,189)$-graphs were produced,
and thus $R(3,10) \le 42$.
\end{proof}

\medskip
Theorem~\ref{theorem:r_3_10_42b} improves over the bound
$R(3,10)\le 43$ obtained in 1988 \cite{RK2}.
The correctness tests of our implementations and
the computational effort required for various parts of
the computations are described in Appendix 2.

\bigskip
Geoffrey Exoo \cite{ExP} found almost 300000 $(3,10;39)$-graphs, we
extended this set to more than $4\cdot 10^7$ graphs, and 
very likely there are more of them. The known $(3,10;39)$-graphs
have the number of edges ranging from 161 to 175, hence we
have $151 \le e(3,9,39) \le 161$. We expect that the actual
value is much closer, if not equal, to 161. Despite many attempts
by Exoo, us, and others, no $(3,10;40)$-graphs were constructed.
The computations needed for the upper bound in Theorem 6 were
barely feasible. Consequently, we anticipate that any further
improvement to either of the bounds in
$40 \le R(3,10) \le 42$ will be very difficult.

\eject

\section{Lower Bounds for $e(3,k,n)$ and\\
Upper Bounds for $R(3,k)$, for $k \ge 11$}
\label{section:larger_10}

\medskip
We establish five further new upper bounds on the Ramsey numbers $R(3,k)$
as listed in Theorem~\ref{theorem:new_upper_bounds}. All of the new bounds
improve the results listed in the 2011 revision of the survey
\cite{SRN} by 1. The bound $R(3,16) \le 98$ was also obtained by Backelin,
though it was not published \cite{Back,BackP}.
Note that we don't improve the upper bound on $R(3,12)$.

\medskip
\begin{thm}
\label{theorem:new_upper_bounds}
The following upper bounds hold:\\[0.8ex]
$R(3,11) \le 50$, $R(3,13) \le 68$,
$R(3,14) \le 77$, $R(3,15) \le 87$, and $R(3,16) \le 98$.
\end{thm}

\medskip
\begin{proof}
Each of the new upper bounds $R(3,k)\le n$ can be obtained by
showing that $e(3,k,n)=\infty$. The details of the intermediate
stages of computations for all $k$ are presented in the tables
and comments of the remaining part of this section.
For $k=16$ no data is shown except
some comments in Table~\ref{table:bounds_e_15},
in particular the data in this table
implies $e(3,16,98)=\infty$ by (4).
\end{proof}

\medskip
In the Tables~\ref{table:bounds_e_11}, \ref{table:bounds_e_12}
and \ref{table:bounds_e_13}, for $k=11, 12$ and 13, respectively,
we list several cases in the comments column, where the lower bounds
on $e(3,k,n)$ listed in \cite{Les} (some of them credited to
\cite{Back}) are better than our results. This is the case
for $n$ slightly larger than $13k/4-1$, mostly due to the theorems
claimed in the unpublished manuscript by Backelin \cite{Back, BackP}.
Our lower bounds on $e(3,k,n)$, and implied upper bounds on $R(3,k)$,
do not rely on these results. We have checked that assuming
the results from \cite{Back, BackP, Les} would not imply, using
the methods of this paper, any further improvements on the upper
bounds on $R(3,k)$ for $k \le 16$, but they may for $k \ge 17$.
Hence, if the results in \cite{Back, Les} are published, then
using them jointly with our results may lead to better
upper bounds on $R(3,k)$, at least for some $k \ge 17$.

\eject
\bigskip
\noindent
{\bf Lower bounds for $e(3,11,n)$}

\medskip
\noindent
The exact values of $e(3,11,\le 31)$ are determined by Theorem~\ref{theorem:comulative_small}.
The bounds for $n=32, 33$ marked with a 't' are from Theorem~\ref{theorem:comulative_small}.
The lower bounds on $e(3,11,\ge 32)$
are included in the second column of Table~\ref{table:bounds_e_11}.
They are based on solving
integer constraints (4), using known values
and lower bounds on $e(3,10,n)$ listed in Table~\ref{table:bounds_e_10} in
Section~\ref{section:computations_10}. They are better than those in \cite{Les}
for all $36 \le n \le 50$.

\begin{table}[H]
\begin{center}
\begin{tabular}{c|c|l}
\hline
$n$&$e(3,11,n)\ge$&comments\cr
\hline
32&\ \ 62t&63 \cite{Les}, 63 is exact \cite{Back,BackP}\cr
33&\ \ 68t&69 \cite{Les}, 70 is exact \cite{Back,BackP}\cr
34&\ 75&76 \cite{Les}, 77 is exact \cite{Back,BackP}\cr
35&\ 83&84 \cite{Les}, credit to \cite{Back}\cr
36&\ 92&the same as in \cite{Back,BackP}\cr
37&100&\cr
38&109&\cr
39&117&unique degree sequence solution, 6-regular\cr
40&128&\cr
41&138&\cr
42&149&\cr
43&159&\cr
44&170&\cr
45&182&\cr
46&195&199 required for $R(3,12) \le 58$\cr
47&209&\cr
48&222&unique solution: $n_9=36, n_{10}=12$,\cr
&&215 required for $R(3,12) \le 59$, old bound\cr
49&237&245 maximum\cr
50&$\infty$&hence $R(3,11) \le 50$, new bound, Theorem~\ref{theorem:new_upper_bounds}\cr
51&$\infty$&hence $R(3,11) \le 51$, old bound\cr
\hline
\end{tabular}

\caption{Lower bounds on $e(3,11,n)$, for $n\ge 32$.}
\label{table:bounds_e_11}

\end{center}
\end{table}

\medskip
The maximum number of edges in any $(3,11;49)$-graph is
that of a 10-regular graph, so a proof of $e(3,11,49)>245$
would imply $R(3,11)\le 49$.
Observe that any graph $G_v$ of any 10-regular $(3,11;50)$-graph
must be a $(3,10;39,150)$-graph. Thus, our improvement of
the upper bound on $R(3,11)$ from 51 to 50 is mainly due
to the new lower bound $e(3,10,39) \ge 151$ (together
with not-too-much-off adjacent bounds).

\eject
\bigskip
\noindent
{\bf Lower bounds for $e(3,12,n)$}

\medskip
\noindent
The exact values of $e(3,12,\le 34)$ are determined by Theorem~\ref{theorem:comulative_small}.
The bounds for $35 \le n \le 37$ marked with a 't' are from Theorem~\ref{theorem:comulative_small}.
The lower bounds on $e(3,12,\ge 35)$
are included in the second column of Table~\ref{table:bounds_e_12}.
They are based on solving
integer constraints (4), using known values
and lower bounds on $e(3,11,n)$ given in Table~\ref{table:bounds_e_11}.
They are better than those in \cite{Les}
for all $43 \le n \le 58$.

\medskip
An improvement of the upper bound on $R(3,12)$ obtained by
Lesser \cite{Les} from 60 to 59 is now immediate (it formed a
significant part of her thesis), but a further improvement
from 59 to 58 would require an increase of the lower bound
on $e(3,12,58)$ by 4.

\begin{table}[H]
\begin{center}
\begin{tabular}{c|c|l}
\hline
$n$&$e(3,12,n)\ge$&comments\cr
\hline
35&\ \ 67t&68 \cite{Les}, 68 is exact \cite{Back,BackP}\cr
36&\ \ 73t&74 \cite{Les}, 75 is exact \cite{Back,BackP}\cr
37&\ \ 79t&81 \cite{Les}, 82 is exact \cite{Back,BackP}\cr
38&\ 86&88 \cite{Les}, 89 \cite{Back}\cr
39&\ 93&95 \cite{Les}, 96 \cite{Back}\cr
40&100&102 \cite{Les}, 103 \cite{Back}\cr
41&109&111 \cite{Les}\cr
42&119&119 \cite{Les}, 120 in \cite{Back}\cr
43&128&the same as in \cite{Back}\cr
44&138&\cr
45&148&\cr
46&158&\cr
47&167&168, proof based on Table 7 \cite{BackP}\cr
48&179&180, proof based on Table 7 \cite{BackP}\cr
49&191&\cr
50&203&\cr
51&216&\cr
52&229&\cr
53&241&\cr
54&255&259 required for $R(3,13) \le 67$\cr
55&269&265 required for $R(3,13) \le 68$, Theorem~\ref{theorem:new_upper_bounds}\cr
56&283&\cr
57&299&\cr
58&316&319 maximum\cr
59&$\infty$&hence $R(3,12) \le 59$, old bound\cr
\hline
\end{tabular}

\caption{Lower bounds on $e(3,12,n)$, for $n\ge 35$.}
\label{table:bounds_e_12}

\end{center}
\end{table}

\eject
\noindent
{\bf Lower bounds for $e(3,13,n)$}

\medskip
\noindent
The exact values of $e(3,13,\le 39)$ are determined by Theorem~\ref{theorem:comulative_small}.
The bound for $n=40$ is from Theorem~\ref{theorem:comulative_small}.
The lower bounds on $e(3,13,\ge 40)$
are included in the second column of Table~\ref{table:bounds_e_13}.
They are based on solving
integer constraints (4), using lower bounds
on $e(3,12,n)$ listed in Table~\ref{table:bounds_e_12}.
They are better than those in \cite{Les}
for all $51 \le n \le 68$.

\begin{table}[H]
\begin{center}
\begin{tabular}{c|c|l}
\hline
$n$&$e(3,13,n)\ge$&comments\cr
\hline
40&\ \ 84t&\ \ 86 \cite{Les}, 87 is exact \cite{BackP}\cr
41&\ 91&\ \ 93 \cite{Les}, 94 is exact \cite{Back}\cr
42&\ 97&100 \cite{Les}, 101 \cite{Back}\cr
43&104&107 \cite{Les}, 108 \cite{Back}\cr
44&112&114 \cite{Les}, 115 \cite{Back}\cr
45&120&122 \cite{Les}, 123 \cite{Back}\cr
46&128&130 \cite{Les}, 132 \cite{Back}\cr
47&136&139 \cite{Les}, 140 \cite{Back}\cr
48&146&148 \cite{Les}\cr
49&157&158 \cite{Les}\cr
50&167&167 \cite{Les}, 168 \cite{Back}\cr
51&177&178 \cite{Back}\cr
52&189&\cr
53&200&\cr
54&212&\cr
55&223&\cr
56&234&\cr
57&247&\cr
58&260&\cr
59&275&\cr
60&289&\cr
61&303&\cr
62&319&326 required for $R(3,14) \le 76$\cr
63&334&\cr
64&350&345 required for $R(3,14) \le 77$, Theorem~\ref{theorem:new_upper_bounds}\cr
65&365&\cr
66&381&\cr
67&398&402 maximum\cr
68&$\infty$&hence $R(3,13) \le 68$, new bound\cr
69&$\infty$&hence $R(3,13) \le 69$, old bound\cr
\hline
\end{tabular}

\caption{Lower bounds on $e(3,13,n)$, for $n\ge 40$.}
\label{table:bounds_e_13}

\end{center}
\end{table}

\eject
\noindent
{\bf Lower bounds for $e(3,14,n)$}

\medskip
\noindent
The exact values of $e(3,14,\le 41)$
are determined by
Theorem~\ref{theorem:comulative_small}.
Only lower bounds on $e(3,14,\ge 66)$
are included in the second column of Table~\ref{table:bounds_e_14},
since these are relevant for our further analysis of
$R(3,15)$ and $R(3,16)$.
They are based on solving
integer constraints (4), using lower bounds
on $e(3,13,n)$ listed in Table~\ref{table:bounds_e_13}.
They are better than those in \cite{Les} for all $66 \le n \le 77$.

\medskip
\begin{table}[H]
\begin{center}
\begin{tabular}{c|c|l}
\hline
$n$&$e(3,14,n)\ge$&comments\cr
\hline
66&321&\cr
67&334&335, proof based on Table 9 \cite{BackP}\cr
68&350&\cr
69&365&\cr
70&381&\cr
71&398&407 required for $R(3,15) \le 86$\cr
72&415&414 required for $R(3,15) \le 87$, Theorem~\ref{theorem:new_upper_bounds}\cr
73&432&\cr
74&449&\cr
75&468&\cr
76&486&494 maximum\cr
77&$\infty$&hence $R(3,14) \le 77$, new bound\cr
78&$\infty$&hence $R(3,14) \le 78$, old bound\cr
\hline
\end{tabular}

\caption{Lower bounds on $e(3,14,n)$, for $n\ge 66$.}
\label{table:bounds_e_14}

\end{center}
\end{table}

\eject
\medskip
\noindent
{\bf Lower bounds for $e(3,15,n)$}

\medskip
\noindent
The exact values of $e(3,15,\le 44)$ are determined by Theorem~\ref{theorem:comulative_small}.
Only lower bounds on $e(3,15,\ge 81)$
are included in the second column of Table~\ref{table:bounds_e_15},
since these are relevant for further analysis of $R(3,16)$.
They are based on solving
integer constraints (4), using lower bounds
on $e(3,14,n)$ listed in Table~\ref{table:bounds_e_14}.
They are better than those in \cite{Les} for all $81 \le n \le 87$.

\medskip
\begin{table}[H]
\begin{center}
\begin{tabular}{c|c|l}
\hline
$n$&$e(3,15,n)\ge$&comments\cr
\hline
81&497&498, proof based on Table 10 \cite{BackP}\cr
82&515&518 required for $R(3,16) \le 97$\cr
  &   &511 required for $R(3,16) \le 98$, Theorem~\ref{theorem:new_upper_bounds}\cr
83&533&\cr
84&552&\cr
85&572&\cr
86&592&602 maximum\cr
87&$\infty$&hence $R(3,15) \le 87$, new bound\cr
88&$\infty$&hence $R(3,15) \le 88$, old bound\cr
\hline
\end{tabular}

\caption{Lower bounds on $e(3,15,n)$, for $n\ge 81$.}
\label{table:bounds_e_15}

\end{center}
\end{table}

\subsection*{Acknowledgements}

We are very grateful to J\"{o}rgen Backelin for careful reading the draft,
numerous comments improving the presentation, his patient guidance through
the maze of lower bounds on $e(3,k,n)$, and the last minute improvements to
some of the bounds (see the comments column of Tables 5, 8, 10 and 11 with
credit to \cite{BackP}).

This work was carried out using the Stevin Supercomputer
Infrastructure at Ghent University. Jan Goedgebeur is
supported by a Ph.D. grant from the Research Foundation
of Flanders (FWO).

\eject
\section*{Appendix 1: Graph Counts}
\label{section:appendix1}

\smallskip
Tables 12--15 below contain all known exact counts
of $(3,k;n,e)$-graphs for specified $n$,
for $k=7,8,9$ and $10$, respectively.
All graph counts were obtained by the algorithms described in
Section~\ref{section:algorithms}. Empty entries indicate 0.
In all cases, the maximum number of edges
is bounded by $\Delta(G)n/2 \le (k-1)n/2$.
All $(3,\le 9;n,\le e(3,k,n)+1)$-graphs
which were constructed by our programs can be obtained from
the \textit{House of Graphs}~\cite{HOG} by searching for the
keywords ``minimal ramsey graph'' or from~\cite{hog-static}.

\bigskip
\begin{table}[H]
\begin{center}
{\scriptsize
\begin{tabular}{| c | ccccccc |}
\hline 
edges & \multicolumn{7}{c}{number of vertices $n$} \vline\\
  $e$ & 16 & 17 & 18 & 19 & 20 & 21 & 22\\
\hline 
20  &  2  &    &    &    &    &    &\\
21  &  15  &    &    &    &    &    &\\
22  &  201  &    &    &    &    &    &\\
23  &  2965  &    &    &    &    &    &\\
24  &  43331  &    &    &    &    &    &\\
25  &  498927  &  2  &    &    &    &    &\\
26  &  4054993  &  30  &    &    &    &    &\\
27  &  ?  &  642  &    &    &    &    &\\
28  &  ?  &  13334  &    &    &    &    &\\
29  &  ?  &  234279  &    &    &    &    &\\
30  &  ?  &  2883293  &  1  &    &    &    &\\
31  &  ?  &  ?  &  15  &    &    &    &\\
32  &  ?  &  ?  &  382  &    &    &    &\\
33  &  ?  &  ?  &  8652  &    &    &    &\\
34  &  ?  &  ?  &  160573  &    &    &    &\\
35  &  ?  &  ?  &  2216896  &    &    &    &\\
36  &  ?  &  ?  &  ?  &    &    &    &\\
37  &  ?  &  ?  &  ?  &  11  &    &    &\\
38  &  ?  &  ?  &  ?  &  417  &    &    &\\
39  &  ?  &  ?  &  ?  &  10447  &    &    &\\
40  &  ?  &  ?  &  ?  &  172534  &    &    &\\
41  &  ?  &  ?  &  ?  &  1990118  &    &    &\\
42-43  &  ?  &  ?  &  ?  &  ?  &    &    &\\
44  &  ?  &  ?  &  ?  &  ?  &  15  &    &\\
45  &  ?  &  ?  &  ?  &  ?  &  479  &    &\\
46  &  ?  &  ?  &  ?  &  ?  &  10119  &    &\\
47  &  ?  &  ?  &  ?  &  ?  &  132965  &    &\\
48  &  ?  &  ?  &  ?  &  ?  &  1090842  &    &\\
49-50  &    &  ?  &  ?  &  ?  &  ?  &    &\\
51  &    &  ?  &  ?  &  ?  &  ?  &  4  &\\
52  &    &    &  ?  &  ?  &  ?  &  70  &\\
53  &    &    &  ?  &  ?  &  ?  &  717  &\\
54  &    &    &  ?  &  ?  &  ?  &  5167  &\\
55  &    &    &    &  ?  &  ?  &  27289  &\\
56  &    &    &    &  ?  &  ?  &  97249  &\\
57  &    &    &    &  ?  &  ?  &  219623  &\\
58  &    &    &    &    &  ?  &  307464  &\\
59  &    &    &    &    &  ?  &  267374  &\\
60  &    &    &    &    &  ?  &  142741  &  1\\
61  &    &    &    &    &    &  43923  &  6\\
62  &    &    &    &    &    &  6484  &  30\\
63  &    &    &    &    &    &  331  &  60\\
64  &    &    &    &    &    &    &  59\\
65  &    &    &    &    &    &    &  25\\
66  &    &    &    &    &    &    &  10\\
\hline
\end{tabular}
}

\caption{Number of $(3,7;n,e)$-graphs, for $n \ge 16$.}

\label{table:graph_counts_K7}
\end{center}
\end{table}

\begin{table}[H]
\begin{center}
{\scriptsize
\begin{tabular}{| c | ccccccccc |}
\hline 
edges & \multicolumn{9}{c}{number of vertices $n$} \vline\\
  $e$ & 19 & 20 & 21 & 22 & 23 & 24 & 25 & 26 & 27\\
\hline 
25  &  2  &    &    &    &    &    &    &    &    \\
26  &  37  &    &    &    &    &    &    &    &    \\
27  &  763  &    &    &    &    &    &    &    &    \\
28  &  16939  &    &    &    &    &    &    &    &    \\
29  &  ?  &    &    &    &    &    &    &    &    \\
30  &  ?  &  3  &    &    &    &    &    &    &    \\
31  &  ?  &  60  &    &    &    &    &    &    &    \\
32  &  ?  &  1980  &    &    &    &    &    &    &    \\
33  &  ?  &  58649  &    &    &    &    &    &    &    \\
34  &  ?  &  1594047  &    &    &    &    &    &    &    \\
35  &  ?  &  ?  &  1  &    &    &    &    &    &    \\
36  &  ?  &  ?  &  20  &    &    &    &    &    &    \\
37  &  ?  &  ?  &  950  &    &    &    &    &    &    \\
38  &  ?  &  ?  &  35797  &    &    &    &    &    &    \\
39  &  ?  &  ?  &  1079565  &    &    &    &    &    &    \\
40-41  &  ?  &  ?  &  ?  &    &    &    &    &    &    \\
42  &  ?  &  ?  &  ?  &  21  &    &    &    &    &    \\
43  &  ?  &  ?  &  ?  &  1521  &    &    &    &    &    \\
44  &  ?  &  ?  &  ?  &  72353  &    &    &    &    &    \\
45  &  ?  &  ?  &  ?  &  2331462  &    &    &    &    &    \\
46-48  &  ?  &  ?  &  ?  &  ?  &    &    &    &    &    \\
49  &  ?  &  ?  &  ?  &  ?  &  102  &    &    &    &    \\
50  &  ?  &  ?  &  ?  &  ?  &  8241  &    &    &    &    \\
51  &  ?  &  ?  &  ?  &  ?  &  356041  &    &    &    &    \\
52  &  ?  &  ?  &  ?  &  ?  &  10326716  &    &    &    &    \\
53-55  &  ?  &  ?  &  ?  &  ?  &  ?  &    &    &    &    \\
56  &  ?  &  ?  &  ?  &  ?  &  ?  &  51  &    &    &    \\
57  &  ?  &  ?  &  ?  &  ?  &  ?  &  3419  &    &    &    \\
58  &  ?  &  ?  &  ?  &  ?  &  ?  &  129347  &    &    &    \\
59  &  ?  &  ?  &  ?  &  ?  &  ?  &  3288695  &    &    &    \\
60  &  ?  &  ?  &  ?  &  ?  &  ?  &  64233886  &    &    &    \\
61-64  &  ?  &  ?  &  ?  &  ?  &  ?  &  ?  &    &    &    \\
65  &  ?  &  ?  &  ?  &  ?  &  ?  &  ?  &  396  &    &    \\
66  &  ?  &  ?  &  ?  &  ?  &  ?  &  ?  &  21493  &    &    \\
67  &    &  ?  &  ?  &  ?  &  ?  &  ?  &  613285  &    &    \\
68  &    &  ?  &  ?  &  ?  &  ?  &  ?  &  11946369  &    &    \\
69-72  &    &  ?  &  ?  &  ?  &  ?  &  ?  &  ?  &    &    \\
73  &    &    &  ?  &  ?  &  ?  &  ?  &  ?  &  62  &    \\
74  &    &    &    &  ?  &  ?  &  ?  &  ?  &  1625  &    \\
75  &    &    &    &  ?  &  ?  &  ?  &  ?  &  23409  &    \\
76  &    &    &    &  ?  &  ?  &  ?  &  ?  &  216151  &    \\
77  &    &    &    &  ?  &  ?  &  ?  &  ?  &  1526296  &    \\
78-84  &    &    &    &    &  ?  &  ?  &  ?  &  ?  &    \\
85  &    &    &    &    &    &    &  ?  &  ?  &  4  \\
86  &    &    &    &    &    &    &  ?  &  ?  &  92  \\
87  &    &    &    &    &    &    &  ?  &  ?  &  1374  \\
88  &    &    &    &    &    &    &    &  ?  &  11915  \\
89  &    &    &    &    &    &    &    &  ?  &  52807  \\
90  &    &    &    &    &    &    &    &  ?  &  122419  \\
91  &    &    &    &    &    &    &    &  ?  &  151308  \\
92  &    &    &    &    &    &    &    &    &  99332  \\
93  &    &    &    &    &    &    &    &    &  33145  \\
94  &    &    &    &    &    &    &    &    &  4746  \\
\hline
\end{tabular}
}

\caption{Number of $(3,8;n,e)$-graphs, for $n \ge 19$.}

\label{table:graph_counts_K8}
\end{center}
\end{table}

\begin{table}[H]
\begin{center}
{\scriptsize
\begin{tabular}{| c | cccccccccccc |}
\hline 
edges & \multicolumn{12}{c}{number of vertices $n$} \vline\\
  $e$ & 24 & 25 & 26 & 27 & 28 & 29 & 30 & 31 & 32 & 33 & 34 & 35\\
\hline 
40  &  2  &    &    &    &    &    &    &    & &&&   \\
41  &  32  &    &    &    &    &    &    &    &  &&&  \\
42  &  2089  &    &    &    &    &    &    &    &  &&&  \\
43  &  115588  &    &    &    &    &    &    &    &  &&&  \\
44-45  &  ?  &    &    &    &    &    &    &    & &&&   \\
46  &  ?  &  1  &    &    &    &    &    &    & &&&   \\
47  &  ?  &  39  &    &    &    &    &    &    & &&&   \\
48  &  ?  &  4113  &    &    &    &    &    &    & &&&   \\
49  &  ?  &  306415  &    &    &    &    &    &    & &&&   \\
50-51  &  ?  &  ?  &    &    &    &    &    &    & &&&   \\
52  &  ?  &  ?  &  1  &    &    &    &    &    &  &&&  \\
53  &  ?  &  ?  &  1  &    &    &    &    &    &  &&&  \\
54  &  ?  &  ?  &  444  &    &    &    &    &    & &&&   \\
55  &  ?  &  ?  &  58484  &    &    &    &    &    & &&&   \\
56-60  &  ?  &  ?  &  ?  &    &    &    &    &    &  &&&  \\
61  &  ?  &  ?  &  ?  &  700  &    &    &    &    &  &&&  \\
62  &  ?  &  ?  &  ?  &  95164  &    &    &    &    & &&&   \\
63  &  ?  &  ?  &  ?  &  6498191  &    &    &    &    & &&&   \\
64-67  &  ?  &  ?  &  ?  &  ?  &    &    &    &    & &&&   \\
68  &  ?  &  ?  &  ?  &  ?  &  126  &    &    &    &  &&&  \\
69  &  ?  &  ?  &  ?  &  ?  &  17223  &    &    &    &  &&&  \\
70  &  ?  &  ?  &  ?  &  ?  &  1202362  &    &    &    & &&&   \\
71-76  &  ?  &  ?  &  ?  &  ?  &  ?  &    &    &    &  &&&  \\
77  &  ?  &  ?  &  ?  &  ?  &  ?  &  1342  &    &    &  &&&  \\
78  &  ?  &  ?  &  ?  &  ?  &  ?  &  156686  &    &    & &&&   \\
79-85  &  ?  &  ?  &  ?  &  ?  &  ?  &  ?  &    &    &  &&&  \\
86  &  ?  &  ?  &  ?  &  ?  &  ?  &  ?  &  1800  &    &  &&&  \\
87  &  ?  &  ?  &  ?  &  ?  &  ?  &  ?  &  147335  &    &  &&&  \\
88-94  &  ?  &  ?  &  ?  &  ?  &  ?  &  ?  &  ?  &    &  &&&  \\
95  &  ?  &  ?  &  ?  &  ?  &  ?  &  ?  &  ?  &  560  &  &&&  \\
96  &  ?  &  ?  &  ?  &  ?  &  ?  &  ?  &  ?  &  35154  &  &&&  \\
97-103  & &  ?  &  ?  &  ?  &  ?  &  ?  &  ?  &  ?  &  &&&  \\
104  &    &    &  ?  &  ?  &  ?  &  ?  &  ?  &  ?  &  39 &&& \\
105  &    &    &    &  ?  &  ?  &  ?  &  ?  &  ?  &  952 &&& \\
106  &    &    &    &  ?  &  ?  &  ?  &  ?  &  ?  &  18598 &&& \\
107  &    &    &    &  ?  &  ?  &  ?  &  ?  &  ?  &  234681 &&& \\
108  &    &    &    &  ?  &  ?  &  ?  &  ?  &  ?  &  2104151 &&& \\
109-117  &    &    &    &    &  ?  &  ?  &  ?  &  ?  & ? &&&  \\
118  &    &    &    &    &    &    &  ?  &  ?  & ? & 5 &&  \\
119  &    &    &    &    &    &    &  ?  &  ?  & ? & 69 &&  \\
120  &    &    &    &    &    &    &  ?  &  ?  & ? & $\ge 1223$ &&  \\
121  &    &    &    &    &    &    &    &  ?  & ? & $\ge 13081$ &&  \\
122  &    &    &    &    &    &    &    &  ?  & ? & $\ge 90235$ &&  \\
123  &    &    &    &    &    &    &    &  ?  & ? & $\ge 401731$ &&  \\
124  &    &    &    &    &    &    &    &  ?  & ? & $\ge 1188400$ &&  \\
125  &    &    &    &    &    &    &    &    & ? & $\ge 2366474$ &&  \\
126  &    &    &    &    &    &    &    &    & ? & $\ge 3198596$ &&  \\
127  &    &    &    &    &    &    &    &    & ? & $\ge 2915795$ &&  \\
128  &    &    &    &    &    &    &    &    & ? & $\ge 1758241$ &&  \\
129  &    &    &    &    &    &    &   &   &  & $\ge 673600$ & 1 &  \\
130  &    &    &    &    &    &    &   &   &  & $\ge 153676$ & 4 &  \\
131  &    &    &    &    &    &    &   &  & & $\ge 18502$ & $\ge 15$ &  \\
132  &    &    &    &    &    &    &   &  & & $\ge 922$ & $\ge 40$ &  \\
133  &    &    &    &    &    &    &   &  & && $\ge 54$ &  \\
134  &    &    &    &    &    &    &   &  & && $\ge 43$ &  \\
135  &    &    &    &    &    &    &   &  & && $\ge 20$ &  \\
136  &    &    &    &    &    &    &   &  & && $\ge 7$ &  \\
137-139  &    &    &    &    &    &    &   &  & && &  \\
140  &    &    &    &    &    &    &    &    &   &&& 1 \\
\hline
\end{tabular}
}

\caption{Number of $(3,9;n,e)$-graphs, for $n \ge 24$.}

\label{table:graph_counts_K9}
\end{center}
\end{table}

\begin{table}[H]
\begin{center}
{\scriptsize
\begin{tabular}{| c | cccccc |}
\hline 
edges & \multicolumn{6}{c}{number of vertices $n$} \vline\\
  $e$ & 29 & 30 & 31 & 32 & 33 & 34\\
\hline 
58  &  5  &    &    &    &   & \\
59  &  1364  &    &    &    &  &  \\
60--65  &  ?  &    &    &    &  &  \\
66  &  ?  &  5084  &    &    &  &  \\
67  &  ?  &  1048442  &    &    & &   \\
68--72  &  ?  &  ?  &    &    &  &  \\
73  &  ?  &  ?  &  2657  &    & &   \\
74  &  ?  &  ?  &  580667  &    & &   \\
75--80  &  ?  &  ?  &  ?  &    &   & \\
81  &  ?  &  ?  &  ?  &  6592  &  &  \\
82--89  &  ?  &  ?  &  ?  &  ?  &  &  \\
90  &  ?  &  ?  &  ?  &  ?  &  57099 & \\
91--98  &  ?  &  ?  &  ?  &  ?  &  ? & \\
99  &  ?  &  ?  &  ?  &  ?  &  ? & $\ge 1$\\
$\ge 100$  &  ?  &  ?  &  ?  &  ?  &  ? & ?\\
\hline
\end{tabular}
}

\caption{Number of $(3,10;n,e)$-graphs, for $29 \le n \le 34$.}
\label{table:graph_counts_K10}
\end{center}
\end{table}

We showed that $e(3,10,34)\ge 99$ (see Section 4),
a $(3,10;34,99)$-graph was constructed by Backelin \cite{BackP},
and thus $e(3,10,34)=99$.

\bigskip
\eject
\section*{Appendix 2: Testing Implementations}
\label{section:appendix2}

\bigskip
\noindent
{\bf Correctness}

\medskip
Since most results obtained in this paper rely on computations,
it is very important that the correctness of our programs has
been thoroughly verified. Below we list the main tests
and agreements with results produced by more than one
computation.

\begin{itemize}
\item
For every $(3,k)$-graph which was output by our programs,
we verified that it does not contain an independent set
of order $k$ by using an independent program.

\item
For every $(3,k;n,e(3,k,n))$-graph which was generated by our
programs, we verified that dropping any edge creates an independent
set of order $k$.

\item
For various $(3,k;n,\le e)$-graphs we added up to $f$ edges in all
possible ways to obtain $(3,k;n,\le e + f)$-graphs. For the cases
where we already had the complete set of $(3,k;n,\le e + f)$-graphs
we verified that no new $(3,k;n,\le e + f)$-graphs were obtained.
We used this, amongst other cases, to verify that no new $(3,9;24,\le 43)$,
$(3,9;28,\le 70)$, $(3,9;30,\le 87)$
or $(3,10;30,\le 67)$-graphs
were obtained.

\item
For various $(3,k;n,\le e + f)$-graphs we
dropped one edge in all
possible ways and verified that no new $(3,k;n,\le e + f - 1)$-graphs
were obtained.
We used this technique, amongst other cases, to verify that no new
$(3,9;24,\le 42)$, $(3,9;28,\le 69)$,
$(3,9;33,\le 119)$, $(3,9;34,\le 130)$,
$(3,10;30,66)$ or $(3,10;32,81)$-graphs were obtained.

\item
For various sets of $(3,k+1;n,\le e)$-graphs we took each member
$G$ and constructed from it all $G_v$'s. We then
verified that this did not yield any new
$(3,k;n - deg(v) - 1,\le e - Z(v))$-graphs for the cases where
we have all such graphs.
We performed this test, amongst other cases, on the sets of
$(3,9;28,\le 70)$- and $(3,10;31,\le 74)$-graphs.

\item
Various sets of graphs can be obtained by both the minimum degree
extension method and the neighborhood gluing extension method. We
performed both extension methods for various cases (e.g. to obtain
the sets of $(3,9;24,\le 43)$ and $(3,9;25,\le 48)$-graphs).
In each of these cases the results obtained by both methods were
in complete agreement.

\item
The sets of $(3,7;21,\le 55)$, $(3,7;22)$, $(3,8;26,\le 76)$
and $(3,8;27,\le 88)$-graphs were obtained by both the maximal
triangle-free method~\cite{BGSP} and the neighborhood gluing
extension method. The results were in complete agreement.
As these programs are entirely independent and the output sets
are large, we think
that this provides strong evidence of their correctness.

\item
The counts of $(3,7;16,20)$, $(3,7;17,25)$, $(3,7;18,30)$,
$(3,7;19,37)$, $(3,7;20,44)$, $(3,7;21,51)$, and $(3,7;22,e)$
for all $60 \le e \le 66$,
are confirmed by~\cite{RK1}.

\item
The counts of $(3,7;18,31)$, $(3,7;19,38)$, $(3,7;20,45)$ and
$(3,7;21,\le 53)$-graphs are confirmed by~\cite{RK2}.

\item
The counts of $(3,8;19,25)$, $(3,8;20,30)$, $(3,8;21,35)$
and $(3,9;24,40)$-graphs are confirmed by~\cite{RK3}.

\item
The counts of $(3,7;16,21)$, $(3,7;17,26)$, $(3,8;22,42)$
and $(3,9;25,47)$-graphs are confirmed by~\cite{BackP}.
\end{itemize}

\smallskip
Additional implementation correctness tests of
specialized algorithms described in Section 6
were as follows:
\begin{itemize}
\item
The specialized program of Section 6 was used to
extend $(3,8;26,76)$- to $(3,9;35,140)$-graphs and it
produced the unique $(3,9;35,140)$-graph.
\item
We relaxed the conditions to generate all $(3,9;32,108)$-graphs
from Lemma 5 by dropping the requirement that each vertex of
degree 6 has 3 neighbors of degree 7, and enforcing just one
vertex of degree 7 with exactly one neighbor of degree 6.
This yielded 21602 graphs. We verified that each of these graphs
was indeed already included in the set $\mathcal{X}$,
and that $\mathcal{X}$ does not contain any additional such graphs.
\end{itemize}

Since our results are in complete agreement with previous results
and since all our consistency tests passed, we believe that this
is strong evidence for the correctness of our implementations.

\bigskip
\noindent
{\bf Computation Time}

\medskip
The implementations of extension algorithms described
in Sections 3 and 6 are written in C.
Most computations were performed on a cluster
with Intel Xeon L5520 CPU's at 2.27 GHz, on which
a computational effort of one CPU year can be usually
completed in about 8 elapsed hours. The overall
computational effort of this project
is estimated to be about 50 CPU years, which includes
the time used by a variety of programs. The most cpu-intensive
tasks are listed in the following.

The first phase of obtaining $(3,9;32,108)$-graphs
required about 5.5 CPU years. The bottlenecks
of this phase were 
the computations required for extending all $(3,8;24,\le 59)$-graphs 
(which required approximately 3.5 CPU years), and extending the 
$(3,8;25,\le 68)$-graphs (which took more than 2 CPU years).
The second phase of obtaining the special $(3,9;32,108)$-graphs with 
$n_6=8$, $n_7=24$ as in Lemma 5 took about 5.8 CPU years.
The specialized program of Section 6 extended all
$(3,9;32,108)$-graphs to 9-regular $(3,10;42,189)$-graphs
quite fast, in about only 0.25 CPU years.
Performing computations to generate all 
$(3,10;39,\le 150)$-graphs (there are none of these),
which were needed for the bound $R(3,11)\le 50$,
took about 4.8 CPU years.

The CPU time needed to complete the computations
of Section 7 was negligible, however their variety
caused that they were performed during the span
of several weeks.

\end{document}